\documentclass[12pt]{article}

 \normalsize

 \usepackage{amssymb}		
 \usepackage{amstext}
 \usepackage{amsfonts}
 \usepackage{amsmath}		
 \usepackage{amscd}		
 \usepackage{latexsym}		
 \usepackage{epsfig}
 \usepackage{subfigure}
 \usepackage{pstricks,pst-node,pst-text,pst-tree}
 \usepackage{array}

 \numberwithin{equation}{section}
 \newtheorem{theorem}{Theorem}[section]
 \newtheorem{prop}[theorem]{Proposition}
 \newtheorem{cor}[theorem]{Corollary}
 \newtheorem{lemma}{Lemma}
 
 \newcommand{\para}{\medskip \indent}
 \newenvironment{pf}{\paragraph{Proof}}{\par\medskip}
 \newenvironment{pfof}[1]{\paragraph{Proof of #1}}{\par\medskip}
 
 \newcommand{\qed}{\ifhmode\unskip\nobreak\fi\quad\ensuremath\square}

 \newcommand{\one}{\ensuremath{\mathrm{i}}}
 \newcommand{\two}{\ensuremath{\mathrm{ii}}}
 \newcommand{\three}{\ensuremath{\mathrm{iii}}}
 
 \newcommand{\A}{\ensuremath{\mathbb{A}}}
 \newcommand{\C}{\ensuremath{\mathbb{C}}}
 \newcommand{\F}{\ensuremath{\mathbb{F}}}
 
 \newcommand{\PP}{\ensuremath{\mathbb{P}}}
 
 \newcommand{\R}{\ensuremath{\mathbb{R}}}
 \newcommand{\Z}{\ensuremath{\mathbb{Z}}}
 \newcommand{\SL}{\ensuremath{\operatorname{SL}}}

 \newcommand{\AHilb}{\ensuremath{A\text{-}\operatorname{\! Hilb}}}
 \newcommand{\GHilb}{\ensuremath{G\text{-}\operatorname{\! Hilb}}}
 \newcommand{\Hilb}{\operatorname{Hilb}}
 \newcommand{\Hom}{\ensuremath{\operatorname{Hom}}}
 \newcommand{\Spec}{\ensuremath{\operatorname{Spec}}}
 \newcommand{\dP}{\ensuremath{\operatorname{dP_{6}}}}
 
 \newcommand{\ahilb}{\ensuremath{A\operatorname{-Hilb}\:\mathbb{C}^{3}}}
 
 \newcommand{\bij}{\leftrightarrow}

 \newcommand{\diag}{\ensuremath{\operatorname{diag}}}

 \newcommand{\iso}{\cong}
 \newcommand{\hcf}{\operatorname{hcf}}

 \newcommand{\Ahilb}[1]{\ensuremath{A\text{-}\operatorname{\! Hilb}\:\mathbb{C}^{#1}}}
 \newcommand{\Frho}[1]{\ensuremath{\mathcal{F}_{#1}}}
 \newcommand{\Ghilb}[1]{\ensuremath{G\text{-}\operatorname{\! Hilb}\:\mathbb{C}^{#1}}}
 \newcommand{\Span}[1]{\left<#1\right>}

 \newcommand{\al}{\alpha}
 \newcommand{\be}{\beta}
 \newcommand{\ga}{\gamma}
 \newcommand{\De}{\Delta}
 \newcommand{\Si}{\Sigma} 
 \newcommand{\la}{\lambda}
 
 \newcommand{\ep}{\varepsilon}
 
 \newcommand{\ze}{\zeta}
 \newcommand{\Oh}{\mathcal O}
 \newcommand{\sI}{\mathcal I}
 \newcommand{\bm}{\mathbf m}

 \title{How to calculate \ahilb}
 \author{Alastair Craw \and Miles Reid}
 \date{}

\begin{document}
 \maketitle
 
 \begin{abstract}
 Nakamura \cite{N} introduced the $G$-Hilbert scheme $\Ghilb{3}$ for a finite subgroup
$G\subset\SL(3,\C)$, and conjectured that it is a crepant resolution of the
quotient $\C^3/G$. He proved this for a diagonal Abelian group $A$ by
introducing an explicit algorithm that calculates $\Ahilb{3}$. This note
calculates $\Ahilb{3}$ much more simply, in terms of fun with continued
fractions plus regular tesselations by equilateral triangles.
 \end{abstract}

 \section{Statement of the result}
 \subsection{The junior simplex and three Newton polygons}
 \label{ssec:3N}
\begin{figure}[thb]
 \centerline{\mbox{\epsfbox{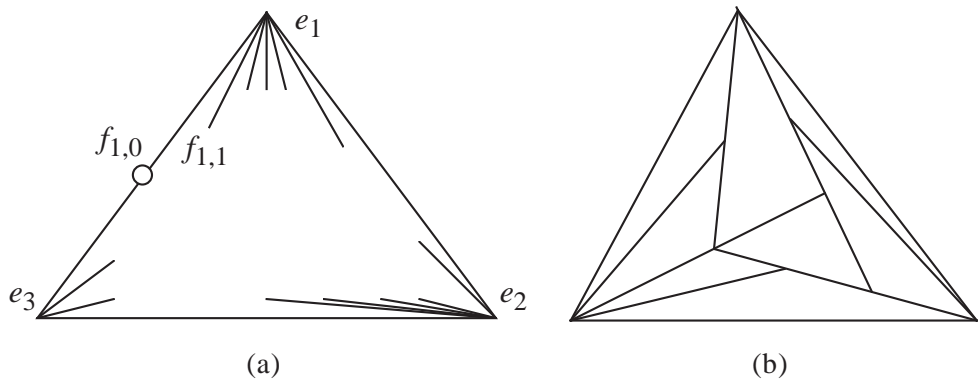}}}
 \caption{(a) Three Newton polygons; (b) subdivision into regular
triangles}
 \label{fig:tri}
 \end{figure}

 Let $A\subset\SL(3,\C)$ be a diagonal subgroup acting on $\C^3$. Write
$L\supset\Z^3$ for the overlattice generated by all the elements of $A$
written in the form $\frac{1}{r}(a_1,a_2,a_3)$. The junior simplex $\De$
(compare \cite{IR}, \cite{R}) has 3 vertexes
 \[
 e_1=(1,0,0),\quad e_2=(0,1,0) \quad\text{and}\quad e_3=(0,0,1).
 \]
Write $\R^2_\De$ for the affine plane spanned by $\De$, and
$\Z^2_\De=L\cap\R^2_\De$ for the corresponding affine lattice. Taking each
$e_i$ in turn as origin, construct the Newton polygons obtained as the
convex hull of the lattice points in $\De \setminus e_{i}$ (see Figure~\ref{fig:tri}.a):
 \begin{equation}
 f_{i,0},f_{i,1},f_{i,2},\dots, f_{i,k_i+1},
 \label{eq:f_ij}
 \end{equation}
 where $f_{i,0}$ is the primitive vector along the side $[e_i,e_{i-1}]$,
and $f_{i,k_i+1}$ that along $[e_i,e_{i+1}]$. (The indices $i,i\pm1$ are
cyclic. Also, since $e_i$ is the origin, the notation $f_{i,j}$
 denotes both the lattice point of $\De$ and the corresponding vector
$e_if_{i,j}$.) The vectors $f_{i,j}$ out of $e_i$ are subject to the
Jung--Hirzebruch continued fraction rule:
 \begin{equation}
 f_{i,j-1}+f_{i,j+1}=a_{i,j} \cdot f_{i,j}
 \quad\text{for $j=1,\dots, k_i$},
 \label{eq:aij}
 \end{equation}
 where $a_{i,j}\ge2$. Here $\frac{r_i}{\al_i}=[a_{i,1},\dots,a_{i,k_i}]$
comes from expressing $\Z^2_\De$ in terms of the cone at $e_i$, writing
 \[
 \Z^2_\De=\Z^2(f_{i,0},f_{i,k_i+1})+\Z\cdot f_{i,1} = \Z^{2} + \Z\cdot
 \textstyle{\frac{1}{r_i}}(\al_i,1),
 \]
 with $\al_i<r$ and coprime to $r$. Write
$L_{ij}$ for the line out of
$e_i$ extending or equal to the initial segment $[e_i,f_{ij}]$ ({\em line}
is {\em line segment}\/ throughout). The resulting fan at $e_i$
 corresponds to the Jung--Hirzebruch resolution of the surface singularity
$\C^2_{(x_i=0)}/A$.  The picture so far is the simplex $\De$ with a
 number of lines $L_{ij}$ growing out of each of the 3 vertexes (Figure~\ref{fig:tri}.a). 

 \subsection{Regular triangles}\label{ssec:reg}
 Write $\Z^2$ for the group of translations of the affine lattice
$\Z^2_\De$. A {\em regular triple} is a set of three vectors
$v_1,v_2,v_3\in\Z^2$, any two of which form a basis of $\Z^2$, and such
that $\pm v_1\pm v_2\pm v_3=0$. (The standard regular triple is
$\pm(1,0),\pm(0,1),\pm(1,1)$; it appears all over elementary toric
geometry, for example, as the fan of $\PP^2$ or the blowup of $\A^{2}$.)
We are only concerned with regular triples among the vectors $f_{i,j}$
introduced in \ref{ssec:3N}.

As usual, a {\em lattice triangle} $T$ is a triangle $T\subset\R_\De^2$
with vertexes in $\Z^2_\De$. We say that $T$ is a {\em regular triangle} if
each of its sides is a line $L_{ij}$ extending some $[e_i,f_{i,j}]$ and
the 3 primitive vectors $v_1,v_2,v_3\in\Z^2$ pointing along its sides form
a regular triple.

\begin{figure}[thb]
 \centerline{\mbox{\epsfbox{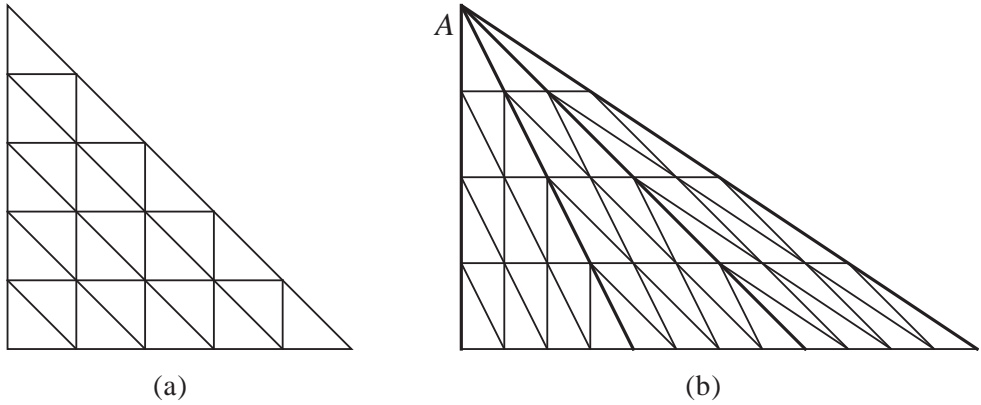}}}
 \caption{(a) A 5-regular triangle; (b) a $(4,12)$-semiregular
triangle (see \ref{sssec:sreg})}
 \label{fig:reg}
 \end{figure}

It is easy to see that a regular triangle $T$ is affine equivalent to the
triangle with vertexes $(0,0),(r,0),(0,r)$ for some $r\ge1$, called the
{\em side} of $T$. Its {\em regular tesselation} is that shown in
Figure~\ref{fig:reg}.a: a regular triangle of side $r$ subdivides into
$r^2$ basic triangles with sides parallel to $v_1,v_2,v_3$.

 A regular triangle is the thing you get as the junior simplex for the
group
 \[
 A=\Z/r\oplus\Z/r=
 \Span{\frac{1}{r}(1,-1,0), \frac{1}{r}(0,1,-1), \frac{1}{r}(-1,0,1)}
 \subset \SL(3,\C)
 \]
(the maximal diagonal subgroup of exponent $r$). The tesselation consists
of basic triangles with vertexes in $\De$, so corresponds to a crepant
resolution of the quotient singularity. It is known (see \ref{ssec:btdmb}
below and \cite{R}, Example~2.2) that in this case $\AHilb\C^3$ is the
toric variety associated with its regular tesselation.

 \subsection{The main result}\label{th:main}

 \begin{theorem}\label{th:1}
 The regular triangles partition the junior simplex $\De$.
 \end{theorem}

 Section~\ref{sec:2} gives an easy continued fraction procedure
determining the partition; Figure~\ref{fig:tri}.b illustrates the
rough idea, and worked examples are given in \ref{ssec:exa}
below\footnote{Homework sheets are on the lecturer's website
www.maths.warwick.ac.uk/$\!\scriptstyle\sim$miles.} (see
Figures~\ref{fig:SmallExa}--\ref{fig:30.25.2.3}).

 \begin{theorem}\label{th:2}
Let \(\Sigma\) denote the toric fan determined by the regular
tesselation (see \ref{ssec:reg}) of all regular triangles in the
junior simplex \(\Delta\). The associated toric variety \(Y_{\Sigma}\) is Nakamura's \(A\)-Hilbert scheme \ahilb.
 \end{theorem}

 \begin{cor}[Nakamura]
 $\AHilb\C^3\to\C^3/A$ is a crepant resolution.
 \end{cor}
 
 \begin{cor}
 \label{cor:excepsurfaces}
 Every compact exceptional surface in \ahilb\ is either \(\mathbb{P}^{2}\),  a scroll \(\F_{n}\) or a scroll blown up in one or two points (including $\dP$,  the del Pezzo surface of degree 6).
 \end{cor}

 \subsection{Thanks} This note is largely a reworking of original ideas of
 Iku Nakamura, and MR had access over several years to his work in progress
 and early drafts of the preprint \cite{N}. MR learned the continued
 fraction tricks here from Jan Stevens (in a quite different context). We
 are grateful to the organisers of two summer schools at Levico in May 1999
 and Lisboa in July 1999 which stimulated our discussion of this material,
 and to Victor Batyrev for the question that we partially answer in
 \ref{sssec:val}.

  \subsection{Recent developments}\label{ssec:recent}

  Since this article first appeared on the e-print server in September
 1999 there has been considerable progress in our understanding of the
 \(G\)-Hilbert scheme.  The most significant development is the work of
 Bridgeland,  King and Reid~\cite{BKR} establishing that
 \(\Ghilb{3}\to \C^{3}/G\) is a crepant resolution for a finite (not
 necessarily Abelian) subgroup \(G\subset \SL(3,\C)\).  In fact \cite{BKR} settles many of the outstanding issues concerning
 \(\Ghilb{3}\);  for instance,  an isomorphism between the K theory of
 \(\Ghilb{3}\) and the representation ring of \(G\) is established,
 and the ``dynamic'' versus ``algebraic'' definition of \(\Ghilb{3}\)
 is settled (see the discussion in Section~\ref{ssec:def} below).

 The explicit calculation of the fan \(\Sigma\) of \(\ahilb\)
  introduced in the current article enabled AC to establish a
  geometric construction of the McKay correspondence.  Indeed,  a certain
  cookery with the Chern classes of the Gonzalez-Sprinberg and Verdier sheaves \(\Frho{\rho}\) (see
  \cite{R} for a discussion) leads to a \(\Z\)-basis of the cohomology \(H^{*}(Y_{\Si},\Z)\) for which the bijection
  \[
  \Big{\{}\mbox{irreducible representations of\ }A\Big{\}}\; \longleftrightarrow\; \mbox{basis of\ } H^*(Y_{\Si},\mathbb{Z})
  \]
  holds,  with \(Y_{\Si} = \ahilb\) (see \cite{C1} for more details).
  Also,  Rebecca Leng's forthcoming Warwick Ph.D.\ thesis \cite{B} extends the explicit calculations in the current
 article to some non-Abelian subgroups of \(\SL(3,\C)\).
 
 Our understanding of the construction of \(\Ghilb{3}\) as a variation
of GIT quotient of \(\C^{3}/G\) has also improved.  Work of King,
Ishii and Craw (summarised in \cite{C2},  Chapter 5) opened the way to
a toric treatment of moduli of representations of the McKay quiver
(also called moduli of \emph{\(G\)-constellations} to stress the link with
\(G\)-clusters).  Initial evidence suggests that these moduli are flops of
\(\Ghilb{3}\):  every flop of \(\Ghilb{3}\) has been constructed in
this way for the quotient of \(\C^{3}\) by the group \(G = \Z/2\times
\Z/2\) (see \ref{ssec:reg}) and for the cyclic quotient singularities \(\frac{1}{6}(1,2,3)\) and \(\frac{1}{11}(1,2,8)\).

 \section{Concatenating continued fractions}\label{sec:2}
 \subsection{Propellor with three blades}\label{ssec:prop}
 The key to Theorem~\ref{th:1} is the observation that easy games with
continued fractions provide all the regular triples $v_1,v_2,v_3$ (see
\ref{ssec:reg}) among the vectors $f_{i,j}$. First translate the three
Newton polygons at $e_1,e_2,e_3$ to a common vertex,
 \begin{figure}[thb]
 \centerline{\mbox{\epsfbox{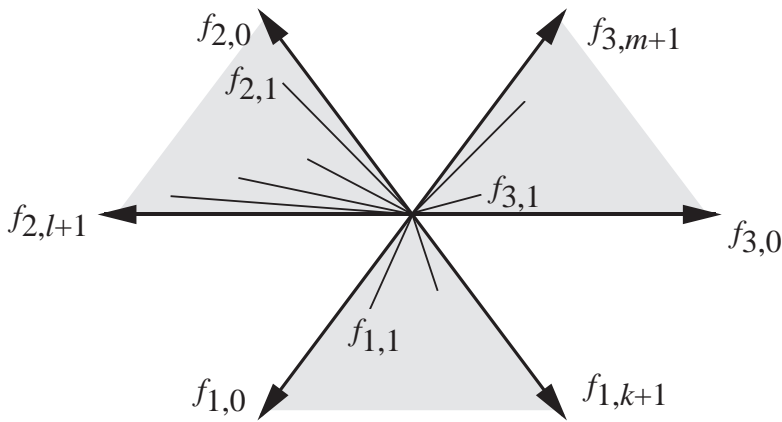}}}
 \caption{``Propellor'' with three ``blades''}
 \label{fig:hex}
 \end{figure}
to get the propellor shape of Figure~\ref{fig:hex}, in which three hexants
(the blades of the propellor) have convex basic subdivisions. The
primitive vectors are read in cyclic order
 \[
 f_{1,0},f_{1,1},\dots,f_{1,k},f_{1,k+1}=-f_{2,0},f_{2,1}, \quad\text{etc.}
 \]
Inverting any blade (that is, multiplying it by $-1$) makes the three
hexants into a basic subdivision of a half-space. Taking plus or minus all
three blades gives a basic subdivision of the plane invariant under $-1$.

 \subsection{Two complementary cones}\label{ssec:2comp}
\begin{figure}[t]
 \centerline{\mbox{\epsfbox{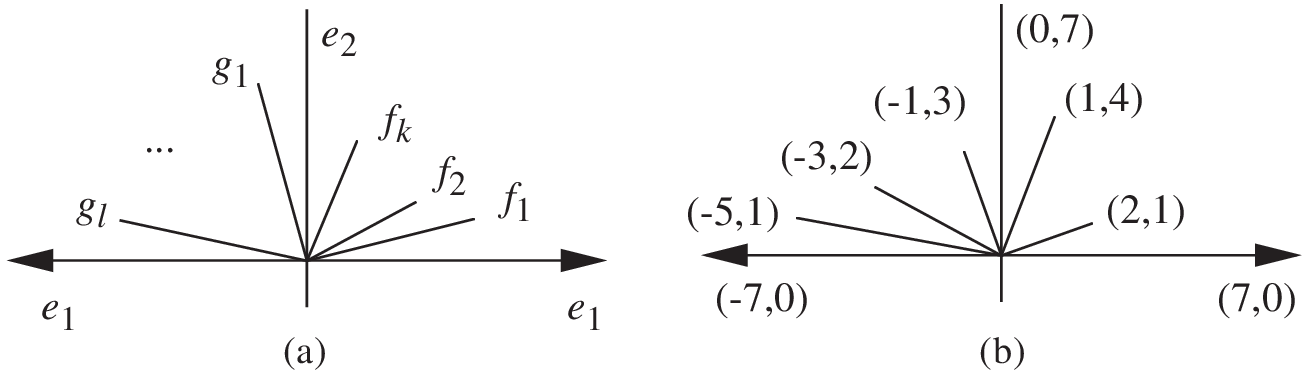}}}
 \caption{Complementary cones $\Span{e_1,e_2}$ and $\Span{e_2,-e_1}$}
 \label{fig:2comp}
 \end{figure}

 This digression on well-known material (see for example \cite{Rie}, \S3,
pp.~220--3) illustrates several points. Let $L$ be a 2-dimensional
lattice, and $e_1,e_2\in L$ primitive vectors spanning a cone in $L_\R$.
Then $\Z^2=\Z\cdot e_1+\Z\cdot e_2\subset L$ is a sublattice with cyclic
quotient $L/\Z^2=\Z/r$; assume for the moment that $r>1$. The {\em
reduced} generator is $f_1=\frac{1}{r}(\al,1)$ with $1\le\al<r$ and
$\al,r$ coprime, so that $L=\Z^2+\Z\cdot\frac{1}{r}(\al,1)$. The continued
fraction expansion $\frac{r}{\al}=[a_1,\dots,a_k]$ with $a_i\ge2$ gives
the convex basic subdivision $\Span{e_1,f_1}$, $\Span{f_i,f_{i+1}}$,
$\Span{f_k,e_2}$ in the first quadrant of Figure~\ref{fig:2comp}.a.

Repeat the same construction for the cone $\Span{e_2,-e_1}$; for this,
write the extra generator $\frac{1}{r}(\al,1)$ as $\frac{1}{r}(\al
e_2,(r-1)(-e_1))$. The reduced normal form is $\frac{1}{r}(1,\be)$ with
$\al\be=(r-1)$ mod $r$, or $\be=1/(r-\al)$ mod $r$. The corresponding
continued fraction $\frac{r}{\be}=[b_1,\dots,b_l]$ gives the basic
subdivision $e_2,g_1,\dots,g_l,-e_1$ in the top left quadrant of
Figure~\ref{fig:2comp}.a. (In the literature, this is usually given as
$\frac{r}{r-\al}=[b_l,\dots,b_1]$, but we want this cyclic order.)

Now the vectors $e_1,f_1,\dots,f_k,e_2,g_1,\dots,g_l,-e_1$ form a basic
subdivision of the upper half-space of $L$. The whole trick is the trivial
observation that this cannot be convex (downwards) everywhere, so that at
$e_2$,
 \begin{equation}
 f_k+g_1=ce_2 \quad\text{with $c\in\Z$ and $0\le c\le 1$.}
 \label{eq:c}
 \end{equation}
For vectors $f_k,g_1$ in the closed upper half-space, $c=0$ is only
possible if $f_k=e_1$ and $g_1=-e_1$. Then $r=1$; this is the ``trivial
case'' with empty continued fractions, at which induction stops. Otherwise,
$f_k+g_1=e_2$. In view of this relation, put a 1 against $e_2$, and
concatenate the two continued fractions as
 \[
 [a_1,a_2,\dots,a_k,1,b_1,\dots,b_l] \quad (=0).
 \]
Because of the relation $e_2=f_k+g_1$, the cone $\Span{f_k,g_1}$ is also
basic. Thus we can delete the vector $e_2$ and still have a basic
subdivision of the upper half-space of $L$. A trivial calculation shows
that in this subdivision, the newly adjacent vectors $f_{k-1},f_k,g_1,g_2$
are related by
 \[
 f_{k-1}+g_1=(a_k-1)f_k \quad\text{and}\quad f_k+g_2=(b_1-1)g_1.
 \]
In other words, in the continued fraction we can replace
 \[
 a_k,1,b_1 \quad\text{by} \quad a_k-1,b_1-1.
 \]
 (The calculation can be seen as the matrix identity
 \[
 \begin{pmatrix}
 0 & 1 \\ -1 & a
 \end{pmatrix}
 \begin{pmatrix}
 0 & 1 \\ -1 & 1
 \end{pmatrix}
 \begin{pmatrix}
 0 & 1 \\ -1 & b
 \end{pmatrix}
 =
 \begin{pmatrix}
 0 & 1 \\ -1 & a-1
 \end{pmatrix}
 \begin{pmatrix}
 0 & 1 \\ -1 & b-1
 \end{pmatrix}.
 \]
The combinatorics is the same as a chain of rational curves on a surface
with self-intersection the negatives of $a_1,a_2,\dots,a_k,
1,b_1,\dots,b_l$; deleting $e_2$ corresponds to ``contracting'' a
$-1$-curve.)

Now it must be the case that at least one of $a_k-1,b_1-1$ is again $1$.
Else the chain of vectors $e_1,f_1,\dots,f_k,g_1,\dots,g_l,-e_1$ is
convex, which is absurd. If say $a_k=2$ then consider the new cone
$\Span{e_1,f_k}$.

Figure~\ref{fig:2comp}.b shows the example $\frac{1}{7}(1,2)$, where we
get
 \begin{equation}
 [4,\underline{2,1,3},2,2] \to [\underline{4,1,2},2,2]
 \to [\underline{3,1,2},2] \to [\underline{2,1,2}] \to [1,1].
 \label{eq:7/4}
 \end{equation}
The steps express $(0,7),(1,4),(-1,3),(-3,2)$ as the sum of two
neighbours. The end $[1,1]$ describes the relations
 \[
 (2,1)=(7,0)+(-5,1)\quad\text{and}\quad (-5,1)=(2,1)+(-7,0)
 \]
among the final four vectors (this counts as one regular triple because we
identify $\pm v$).

\subsection{Remarks}\label{ssec:rem}
 \begin{enumerate}
 \item In the trivial case $r=1$ we have $c=0$ in (\ref{eq:c}). There is
always a 1 to contract. You always end up with $[1,1]=0$.
 \item The regular triples $v_1,v_2,v_3$ among
$e_1,f_1,\dots,e_2,g_1,\dots,-e_1$ correspond one-to-one with the 1's that
occur during the chain of contractions, as we saw in
Figure~\ref{fig:2comp}.b.
 \item The order the vectors are contracted and the regular triples among
them is determined in the course of an induction; but they might be tricky
to decide a priori without running the algorithm.
 \item The continued fractions keep track of successive change of basis
between adjacent basic cones. Following $(e_1,f_1)$, $(f_1,f_2)$, etc.\
all the way around to $(g_l,-e_1)$, and on cyclically to $(-e_1,-f_1)$
gives
 \begin{multline*}
 \begin{pmatrix}-1&0\\0&-1\end{pmatrix}
 =
 \begin{pmatrix}0&1\\-1&a_1\end{pmatrix} \cdots
 \begin{pmatrix}0&1\\-1&a_k\end{pmatrix}
 \begin{pmatrix}0&1\\-1&1\end{pmatrix} \times \\ \times
 \begin{pmatrix}0&1\\-1&b_1\end{pmatrix} \cdots
 \begin{pmatrix}0&1\\-1&b_l\end{pmatrix}
 \begin{pmatrix}0&1\\-1&1\end{pmatrix}.
 \end{multline*}
In what follows, we consider continued fractions concatenated in this cyclic
way. Then $[1,1,1]$ stands for
$\left(\begin{smallmatrix}0&1\\-1&1\end{smallmatrix}\right)^3=
\left(\begin{smallmatrix}-1&0\\0&-1\end{smallmatrix}\right)$, which makes
sense of the number $[1,1,1]=1-\frac{1}{0}=\infty$.
 \end{enumerate}

 \subsection{Long side}\label{ssec:long}
 To concatenate the three continued fractions arising from the propellor
of Figure~\ref{fig:hex} as a cyclic continued fraction, we study the
change of basis from the last basis $f_{1,k},f_{1,k+1}$ of the $e_1$
hexant to the first basis $f_{2,0},f_{2,1}$ of the $e_2$ hexant. Clearly
$f_{2,0}=-f_{1,k+1}$, and we claim there is a relation
 \begin{equation}
 f_{2,1}-f_{1,k}=cf_{2,0} \quad\text{with $c\ge1$}.
 \label{eq:long}
 \end{equation}
Indeed, $-f_{1,k},f_{2,0}$ and $f_{2,0},f_{2,1}$ are two oriented bases
(the usual argument).

\begin{figure}[thb]
 \centerline{\mbox{\epsfbox{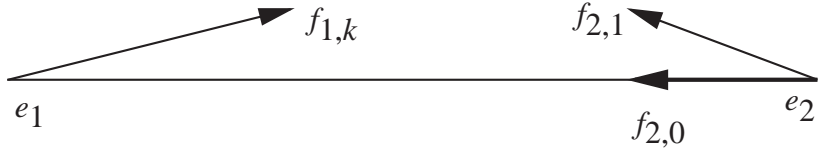}}}
 \caption{A long side of $\De$: $f_{2,1}-f_{1,k}=cf_{2,0}$ with $c\ge2$}
 \label{fig:long}
 \end{figure}

We define the side $e_ie_{i+1}$ of the simplex $\De$ to be a {\em long
side} if $c\ge2$. See Figure~\ref{fig:long}.
A long side $e_1e_2$ is obviously not a primitive vector, so never occurs
for ``coprime'' groups. The presence of a long side is a significant
dichotomy in the construction (see Remark~\ref{ssec:uni}.2).

 \begin{lemma} $\De$ has at most one long side.
 \end{lemma}

If $e_1e_2$ and $e_1e_3$ (say) are both long sides, the basic subdivision
of the upper half-space obtained by inverting the bottom blade of the
propellor in Figure~\ref{fig:hex} would be convex at each ray; this is a
contradiction, as usual. \qed

 \subsection{Concatenating three continued fractions}\label{ssec:conc}
 Suppose that $e_1e_3$ and $e_2e_3$ are not long sides, and that $e_1e_2$
has $c\ge1$ in (\ref{eq:long}). Consider the cyclic continued fraction:
 \begin{equation}
 [1,a_{1,1},\dots,a_{1,k_1},\underline{c},a_{2,1},\dots,a_{2,k_2},
 1,a_{3,1},\dots,a_{3,k_3}].
 \label{eq:cyc}
 \end{equation}
 As above, the meaning of this is the successive change of bases
anti\-clockwise around the figure, from $f_{1,0},f_{1,1}$ to
$f_{1,1},f_{1,2}$ to $f_{1,k},f_{1,k+1}$, then inverting to
$-f_{1,k},f_{1,k+1}=f_{2,0}$ etc., and on to $-f_{1,0},-f_{1,1}$. For most
purposes, we can afford to be sloppy, and not distinguish between $\pm
f_{ij}$, especially in view of the definition of regular triple in
\ref{ssec:reg}. The continued fraction (or any cyclic permutation of it)
evaluates to $\infty=1-\frac{1}{0}$, as explained in
Remark~\ref{ssec:rem}.4.

 \subsection{Examples}\label{ssec:exa}

\begin{figure}[thb]
 \centerline{\mbox{\epsfbox{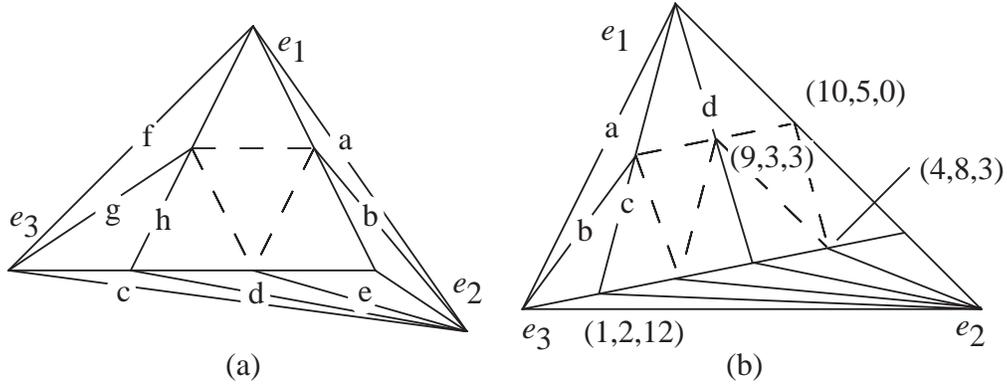}}}
 \caption{Deconstructing (a) $\frac{1}{11}(1,2,8)$ and (b)
$\frac{1}{15}(1,2,12)$: at each step, delete a regular triangle with side
the condemned vector}
 \label{fig:SmallExa}
 \end{figure}

 \paragraph{An example with no long side: $\frac{1}{11}(1,2,8)$} The three
continued fractions (see Figure~\ref{fig:SmallExa}.a) are
 \begin{align*}
 &\text{at $e_1$: \quad $\frac{11}{4}=[3,4]$ \quad
 (because $\frac{1}{11}(2,8)=\frac{1}{11}(1,4)$),} \\
 &\text{at $e_2$: \quad $\frac{11}{7}=[2,3,2,2]$ \quad
 (because $\frac{1}{11}(8,1)=\frac{1}{11}(1,7)$),}\\
 &\text{at $e_3$: \quad $\frac{11}{2}=[6,2]$}.
 \end{align*}
Since the group is coprime, there is
no long side, and these concatenate to
 \begin{equation}
 [1,3,4,1,2,3,2,2,1,6,2] \quad (=\infty).
 \label{eq:128}
 \end{equation}
The contraction rule $\underline{a,1,b}\to a-1,b-1$ is as in
\ref{ssec:2comp}. After any number of contractions, a 1 means a regular
triple $v_1,v_2,v_3$ among the $f_{i,j}$.

Each 1 in (\ref{eq:128}) corresponds to one of the sides $e_3e_1$,
$e_1e_2$ and $e_2e_3$. A chain of contractions with only one 1 allowed to
eat its neighbours corresponds to deleting regular triangles along that
side (see Figure~\ref{fig:SmallExa}.a): contractions along different sides
``commute'', in the sense that they can be done independently of one
another. Thus starting afresh from $[1,3,4,1,2,3,2,2,1,6,2]$ each time
(and numbering the steps as in Figure~\ref{fig:SmallExa}.a), we can do
 \begin{align*}
 & \begin{array}{rrcl}
\text{Step a} & f_{2,0}=f_{2,1}-f_{1,2}: &\to & [1,3,3,1,3,2,2,1,6,2] \\
\text{Step b} & f_{2,1}=f_{2,2}-f_{1,2}: &\to & [1,3,2,2,2,2,1,6,2] \\
 \end{array}\\[6pt]
 \text{or} \quad & \begin{array}[t]{rrcl}
\text{Step c} & f_{2,5}=f_{2,4}-f_{3,1}: &\to & [1,3,4,1,2,3,2,1,5,2] \\
\text{Step d} & f_{2,4}=f_{2,3}-f_{3,1}: &\to & [1,3,4,1,2,3,1,4,2] \\
\text{Step e} & f_{2,3}=f_{2,2}-f_{3,1}: &\to & [1,3,4,1,2,2,3,2] \\
 \end{array}\\[6pt]
 \text{or} \quad & \begin{array}[t]{rrcl}
\text{Step f} & f_{1,0}=f_{1,1}-f_{3,2}: &\to & [2,4,1,2,3,2,2,1,6,1] \\
\text{Step g} & f_{3,2}=f_{3,1}-f_{1,1}: &\to & [1,4,1,2,3,2,2,1,5] \\
\text{Step h} & f_{1,1}=f_{1,2}-f_{3,1}: &\to & [3,1,2,3,2,2,1,4] \\
 \end{array}
 \end{align*}
Carrying out all of these in this order finally gives $[1,1,1]$, which
corresponds to the regular triple $f_{1,2}+f_{2,2}+f_{3,1}=0$. (There is
no uniqueness here, but this is obviously a sensible choice; this
end-point is a meeting of champions as in Remark~\ref{ssec:uni}.2.)

 \paragraph{Example of a long side: $\frac{1}{15}(1,2,12)$}
 Note that $\hcf(15,12)=3$, and the primitive vector along $e_1e_2$ is
$f_{1,3}=-f_{2,0}=(-5,5,0)$ (I omit denominators $\frac{1}{15}$
throughout); see Figure~\ref{fig:SmallExa}.b. Since $f_{1,2}=(-6,3,3)$,
$f_{2,1}=(4,-7,3)$ we see that $f_{2,1}-f_{1,2}=2f_{2,0}$ and $e_1e_2$ is
a long side with $c=2$. In this case, because of the common factor, the
cones at $e_1$ and $e_2$ are $\frac{1}{15}(1,6)\sim\frac{1}{5}(1,2)=[3,2]$
and $\frac{1}{5}(4,1)=[2,2,2,2]$. At $e_3$ we have
$\frac{1}{15}(2,1)=[8,2]$.

Thus the concatenation (\ref{eq:cyc}) is
 \[
 [1,3,2,\underline{2},2,2,2,2,1,8,2].
 \]
A chain of 5 contractions centred around the second 1 corresponds to
deleting the 5 basic triangles along the bottom
Figure~\ref{fig:SmallExa}.b, and reduces the continued fraction to
[1,3,2,1,3,2]. The last of these contractions cuts the long side down to
ordinary size by deleting the bottom right triangle. Alternatively,
starting from the first 1, the 4 steps
 \begin{multline*}
 [1,3,2,2,2,2,2,2,1,8,2] \to
 [2,2,2,2,2,2,2,1,8,1] \\ \to
 [1,2,2,2,2,2,2,1,7] \to
 [1,2,2,2,2,2,1,6] \to
 [1,2,2,2,2,1,5]
 \end{multline*}
deletes the top 4 regular triangles (two of them of side 2) in the order
indicated in Figure~\ref{fig:SmallExa}.b, the last step also cutting
the long side down to size. Doing all of these steps deletes all the
triangles. Note that there are no regular triangles along the long side
$e_1e_2$.

 \subsection{MMPs and regular triples} \label{ssec:MMPs}
 \begin{lemma} For brevity, call a chain of contractions taking a cyclic
continued fraction (\ref{eq:cyc}) down to $[1,1,1]$ an {\em MMP}.
 \begin{enumerate}
 \renewcommand{\labelenumi}{(\roman{enumi})}
 \item Every contraction of\/ $1$ in an MMP corresponds to a regular
triple.
 \item For every regular triple, there is MMP ending at it.
 \item Every regular triple appears in every MMP.
 \end{enumerate}

 \end{lemma}

 \begin{pf} In this proof, view the $\{f_{ij}\}$ as defining a fan of
basic cones invariant under $-1$; we completely ignore the given
``propellor'', and identify $\pm v$.  

 A 1 corresponds to a relation $v_2=v_1+v_3$, which is (i). (ii) is clear: if $v_2=v_1+v_3$ is a
regular triple, then $v_1,v_2,v_3$ and their minuses subdivide $\R^2$
into 6 basic cones. The chain of vectors $f_{ij}$ within any cone is a
nonminimal basic subdivision, so contracts down. 

 We prove (iii): given a regular triple $v_1,v_2,v_3$ and a choice of MMP, suppose
that the first step affecting any of the $v_i$ contracts $v_3$, and choose
signs so that $v_3=v_1+v_2$. Then $v_1,v_2$ span a basic convex cone, and
the original vectors $f_{ij}$ (including $v_3$) form a basic subdivision.
After contracting some of these, the step under consideration contracts
$v_3$, and thus writes it as the sum of two adjacent integral vectors,
which must be in the cones $\Span{v_1,v_3}$ and $\Span{v_2,v_3}$. Since
we're asking for a solution to $(1,1)=(a,b)+(c,d)$ with integers $a>b\ge0$
and $d>c\ge0$, it's clear that the only possible such expression is
$v_3=v_1+v_2$. \qed \end{pf}

\noindent \textbf{Alternative proof of (\three):} Count the number of regular triples and the
number of contractions in an MMP.  It's clear from the MMP algorithm
that each vector \(v_{i}\) appears in precisely \(c_{i}\) regular
triples,  where \(c_{i}\) is the strength of \(v_{i}\).   It follows
that the disjo int union of all regular triangles has \(\sum c_{i}\) edges,  so there are \(\frac{1}{3}\sum c_{i}\) distinct regular triples.  On the other hand,  in a given MMP each contraction reduces the total strength (i.e.\ the sum of the numbers in the continued fraction) by three so there are \(\frac{1}{3}\sum c_{i}\) contractions.  The result follows from the observation that a regular triple cannot correspond to more than one contraction in a given MMP.\qed

 \medskip The lemma says that $\De$ has a unique subdivision into regular triangles,
and any MMP computes it. This completes the proof of Theorem~\ref{th:1}.
\qed

 \subsection{Remarks}\label{ssec:uni}
 Before proceeding to $\GHilb$ and the proof of Theorem~\ref{th:2},
there's still a lot of fun to be derived from regular triples and the
subdivision of Theorem~\ref{th:1}.

 \subsubsection{It's a knock-out!}\label{sssec:KO} The MMP in cyclic
 continued fractions has an entertaining interpretation as a contest
 between the lines $L_{i,j}$ which emanate from the 3 vertices $e_i$.
 The fan \(\Sigma\) of \ahilb\ can be calculated using a
 simple 3-step procedure:

 \begin{enumerate}
 \item Draw lines \(L_{ij}\) emanating from the corners of
 \(\Delta\) (as illustrated in Figure~\ref{fig:tri}.a).  Record the
 strength $a_{ij}$ determined by the Jung--Hirzebruch continued
 fraction rule (\ref{eq:f_ij}) on each line. 
 \item Extend the lines \(L_{ij}\) until they are `defeated' by lines \(L_{kl}\) from \(e_{k}\) (\(i\neq k\)) according to the following rule:  when two or more lines meet at a point,  the line with greater strength extends but its strength decreases by 1 for every rival it defeats.  Lines which meet with equal strength all die.  As a consequence,  strength 2 lines always die.
 \item Step 2 produces the partition of \(\Delta\) into regular triangles of Theorem~\ref{th:1}.  The regular tesselation of the regular triangles gives \(\Sigma\).
 \end{enumerate}

\begin{figure}[!ht]
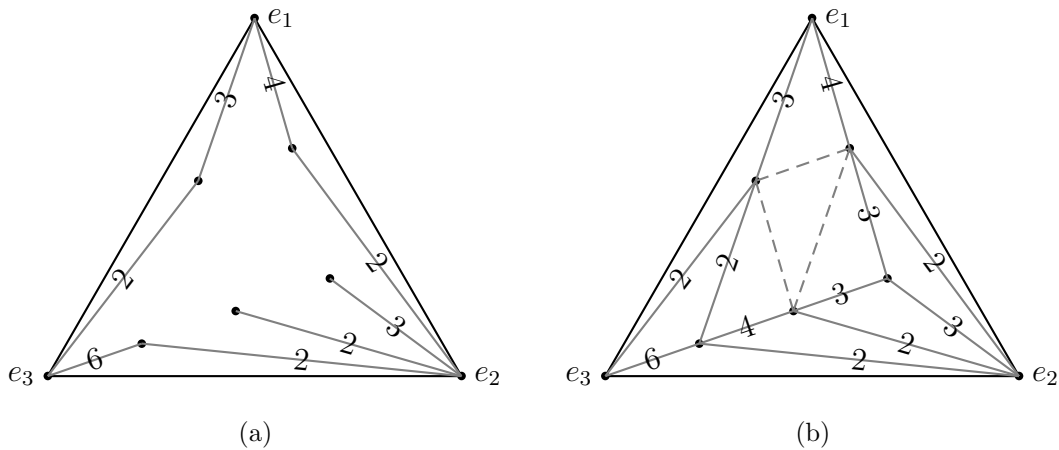

 \centering
 \mbox{\subfigure[]{\input{step1.pspic}\label{fig:11.1.2.8(a)}} \hspace{1.5cm} \subfigure[]{\input{step2.pspic}\label{fig:11.1.2.8(b)}}}
 \caption{(a) Step 1; (b) Step 2 (solid lines) and Step 3 (dotted lines)}
 \label{fig:steps}
 \end{figure}

 \paragraph{Example $\frac{1}{11}(1,2,8)$ revisited:} 
 \label{ex:11.1.2.8}
 Consider the cyclic quotient singularity of type \(\frac{1}{11}(1,2,8)\).  The three continued fractions are
 \[
 \textstyle{\frac{11}{4}=[3,4]\mbox{ at $e_1$;}\quad \frac{11}{7}=[2,3,2,2]\mbox{ at $e_2$;}\quad  \frac{11}{2}=[6,2]\mbox{ at $e_3$.}} 
 \]
 \indent Figure~\ref{fig:11.1.2.8(a)} illustrates the result of Step~1 of the procedure.
 \noindent The solid lines in Figure~\ref{fig:11.1.2.8(b)} show the result of
 Step~2.  For example,  the line from \(e_{1}\) with strength 3
 intersects the line from \(e_{3}\) with strength 2;  the procedure
 says that the line from \(e_{1}\) extends with strength 2 while the
 line from \(e_{3}\) terminates.  The resulting partition of
 \(\Delta\) contains only one regular triangle of side \(r > 1\).  To
 perform Step~3 simply add the dotted lines to Figure~\ref{fig:11.1.2.8(b)}.

\begin{figure}[!ht]
 \centering
 \input{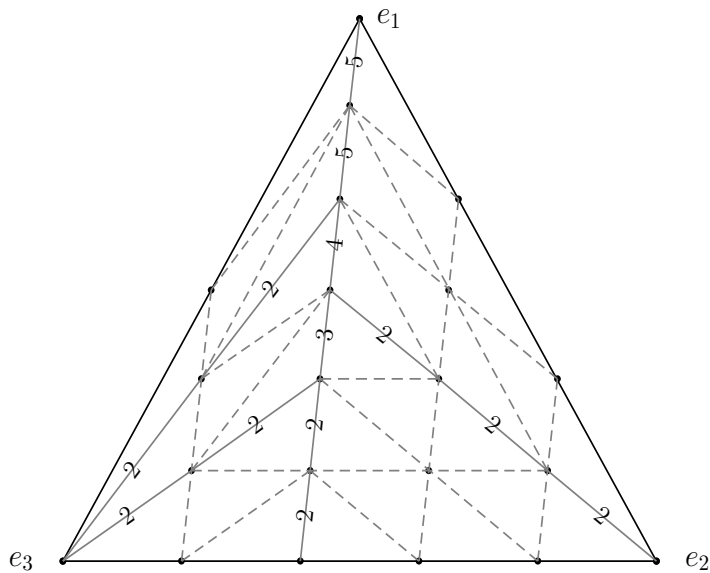}
 \caption{``It's a knock-out!'' for the example $\frac{1}{30}(25,2,3)$}
 \label{fig:30.25.2.3}
 \end{figure}

 \paragraph{Another long sided example: $\frac{1}{30}(25,2,3)$} 
 \label{ex:30.25.2.3}
 Consider the cyclic quotient singularity of type
 \(\frac{1}{30}(25,2,3)\).  Note that \(\hcf(30,25) = 5\) and, because
 of the common factor,  the three continued fractions are
 $\frac{1}{30}(2,3)\sim\frac{1}{5}(1,1) = [5]$ at
 $e_1$,  $\frac{1}{30}(25,2)\sim\frac{1}{2}(1,1) = [2]$ at
 $e_2$ and $\frac{1}{30}(25,2) \sim \frac{1}{3}(2,1) = [2,2]$ at $e_3$.  The solid lines in Figure~\ref{fig:30.25.2.3},  each marked with the
 appropriate strength,  show the partition of the junior simplex of
 \(\frac{1}{30}(25,2,3)\) into regular triangles of side two and
 three.  The dotted lines tesselate the regular triangles.
 
 \para To have some fun, make some extra photocopies of
p.~\pageref{fig:BigExa} to distribute to the class. This is a special
homework sheet doing the example $\frac{1}{101}(1,7,93)$. All the
ideas of the paper can be worked out in detail on it (solutions not provided).

 \subsubsection{Meeting of champions} 
\begin{figure}[thb]
 \centerline{\mbox{\epsfbox{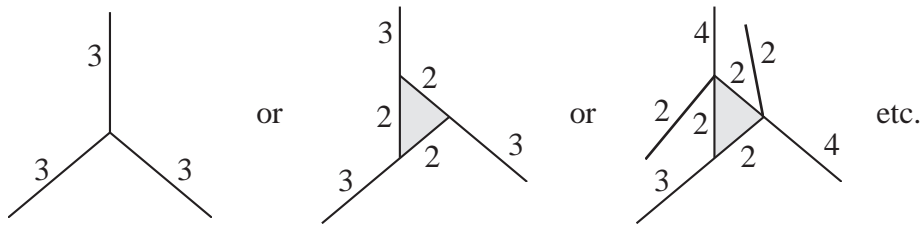}}}
 \caption{Meeting of champions}
 \label{fig:champs}
 \end{figure}

A regular triple is in one of two
possible orientations:
 \begin{description}
 \item[Type 1:] two consecutive vectors in the same closed blade of the
propellor, for example, $f_{1,2}=f_{1,1}+f_{3,1}$ of Figure~\ref{fig:hex};
or
 \item[Type 2:] an interior vector in each blade, for example
$f_{1,2}+f_{2,2}+f_{3,1}=0$.
 \end{description}
If there is a long side $e_1e_2$, it is subdivided by a line from $e_3$,
and Type~2 cannot occur. We claim that if there is no long side, there is
a unique regular triple of Type~2, giving either 3 concurrent vectors or a
cocked hat as in Figure~\ref{fig:champs}; both cases occur (see
Figure~\ref{fig:SmallExa}.a and \cite{R}, Figure~10).
These three are the champions of the knock-out competition, that meet after
eliminating all their less successful rivals.

 \begin{pfof}{claim} Uniqueness is almost obvious from the topology: if it
exists, a meeting of champions divides $\De$ into 4 regions (one possibly
empty), and any other line is confined to one region (it is knocked out by
any champion it meets).

For the existence, the idea is that it is natural to deconstruct $\De$ by
eating in from one side, as we did in the examples of \ref{ssec:exa}. The
cyclic continued fraction (\ref{eq:cyc}) has three 1's, so that each side
of $\De$ takes part in one regular triangle. Choose one side (say
$e_1e_3$) and, preserving the other two, eat as many regular triangles as
we can along $e_1e_3$ (that is, with sides through $e_1$ or $e_3$, as in
Figure~\ref{fig:2reg}.a). Every regular triple of Type~1 is associated
with a well defined side of $\De$, and is eaten in this way starting from
that side. The union of regular triangles along each side forms its {\em
catchment area}\footnote{For example,  in Figure~\ref{fig:30.25.2.3} the three regular triangles of side 2 form the catchment area of \(e_{1}e_{3}\) and the two regular triangles of side 3 form the catchment area of \(e_{1}e_{2}\).  The division into catchment areas determines a `coarse subdivision' of \(\Delta\);  see Craw~\cite{C1}, \S7.1.}.

We now view a MMP as successively deleting dividing lines of the
sub\-division of Figure~\ref{fig:hex}. Eating triangles in the catchment
area of side $e_1e_3$ only deletes lines in the two hexants in the top
right of Figure~\ref{fig:hex}, between $f_{2,0}$ and $f_{3,0}$. Deleting a
line joins two old cones to make a new cone, which is always basic; we
conclude that the two vectors $v,v'$ bounding the catchment area of
$e_1e_3$ form a basis. After this, by assumption, no remaining line in
these two hexants is marked with $1$, so that the cone
$\Span{f_{2,0},f_{3,0}}$ now has its standard Newton polygon subdivision.

If we now complete an MMP anyhow from this position, the same two vectors
$v,v'$ must occur in some regular triple. By what we have said, the
remaining vector must be in the interior of the third hexant. This proves
that a regular triple of Type~2 exists. \qed \end{pfof}

 \subsubsection{Semiregular triangles}\label{sssec:sreg}
 The following definition is not logically part of
Theorems~\ref{th:1}--\ref{th:2}, but it helps to understand complicated
examples: a triangle $T=\De ABC$ (with preferred vertex $A$) is
$(r,cr)$-{\em semi\-regular} if it is equivalent to the triangle with
vertexes $(r,0),(0,0),(0,cr)$. Its {\em semi\-regular tesselation} is
that shown in Figure~\ref{fig:reg}.b. View a $(r,cr)$-semiregular triangle
as made up of $c$ adjacent $r$-regular triangles with vertex at $A$; its
semiregular triangulation is obtained by taking regular triangulations of
each of these. (Note that we work with the affine group of $\Z^2$, so that
each regular triangulation is a perspective view of a tesselation by
equilateral triangles.) If $v_1,v_2,v_3$ are the primitive vectors along
the sides of $T$ (in cyclic order, with $v_1$ the preferred side opposite
$A$), the diagnostic test for semiregularity is that $v_1,v_2$ base
$\Z_\De$ and $cv_1+v_2+v_3=0$. A semiregular triangle relates in the same
way as in \ref{ssec:reg} above to the group $\Z/r\oplus\Z/cr=
\Span{\frac{1}{r}(1,-1,0),\frac{1}{cr}(0,1,-1)}$. The cyclic continued
fraction of a $(r,cr)$-semiregular triangle is $[1,2,2,\dots,2,1,c]$
 with a chain of $c-1$ repeated 2's.

The point of the definition is that it allows you to ignore a string of 2's
in continued fractions. If you calculate a series of examples such as
$\frac{1}{101}(1,k,100-k)$ for $k=2,3,4,5,6$ you'll see that almost all the
area of $\De$ is taken up by semiregular triangles, so this definition is
a convenient way of summarising the information.

In this kind of toric geometry, the following objects correspond: (1) a
string of 2's in a continued fraction; (2) the continued fraction of
$\frac{r}{r-1}$ and the matrix $\left(\begin{smallmatrix} r-1 & r-2 \\ r &
r-1 \end{smallmatrix}\right)$; (3) a row of collinear points in $L$; (4) a
chain of $-2$-curves; (5) an $A_k$ singularity on the relative canonical
model of a surface.

 \subsubsection{Description of $\Si$}
 \label{sssec:val}
 It is not hard to read from the construction of the basic fan $\Si$ that
every (internal) vertex has valency $3,4,5$ or 6, and every (compact)
surface of the resolution is $\PP^2$, a scroll $\F_{n}$, or a once or twice
blown-up scroll including $\dP$ (the del Pezzo surface of degree 6,
 the regular hexagons of \cite{R}).  This provides the foundation for
 an explicit construction of the McKay correspondence for \(\ahilb\) (see \cite{C1}).  The $\dP$ correspond to internal lattice points in
the tesselations of the regular triangles; there are $\binom{r_i-1}{2}$ of
them in each regular triangle of side $r_i$.  Looking at what happens
 in examples,  including quite complicated ones (see the Activity Pack
 on p.~\pageref{fig:BigExa}), seems to indicate other restrictions on $\Si$: for example, a twice blown up scroll usually has a twice blown up fibre
with 3 components of selfintersection $-2,-1,-2$; scrolls $\F_{a}$ or blown
up scrolls only glue into other $\F_{a'}$ with $|a-a'|\le2$. This
question deserves a more systematic study.

 \subsubsection{Inflation and further regular subdivision}\label{ssec:infl}
Note that inflating $\De$ to $n\De$ (or equivalently, replacing $\Z^2_\De$
by $\frac{1}{n}\Z^2_\De$), which corresponds to extending $A$ to
$n^2A=\{g\in\diag\cap\SL(3,\C)\bigm|ng\in A\}$, leaves the continued
fractions at the corners unchanged, so the same picture still gives a
subdivision into regular triangles, with a finer meshed regular
tesselation.


 \section{Regular triangles versus invariant ratios of monomials}
 \label{sec:2reg}
 \subsection{Regular triples and invariant ratios}
 \label{ssec:regtrip}
 The regular triples $v_1,v_2,v_3$ of Section~\ref{sec:2} live in $L$.
Passing to the dual lattice $M$ of invariant monomials is a clever
exercise in elementary coordinate geometry in an affine lattice that plays
a key role in the proof of Theorem~\ref{th:2}.

\begin{figure}[htb]
 \centerline{\mbox{\epsfbox{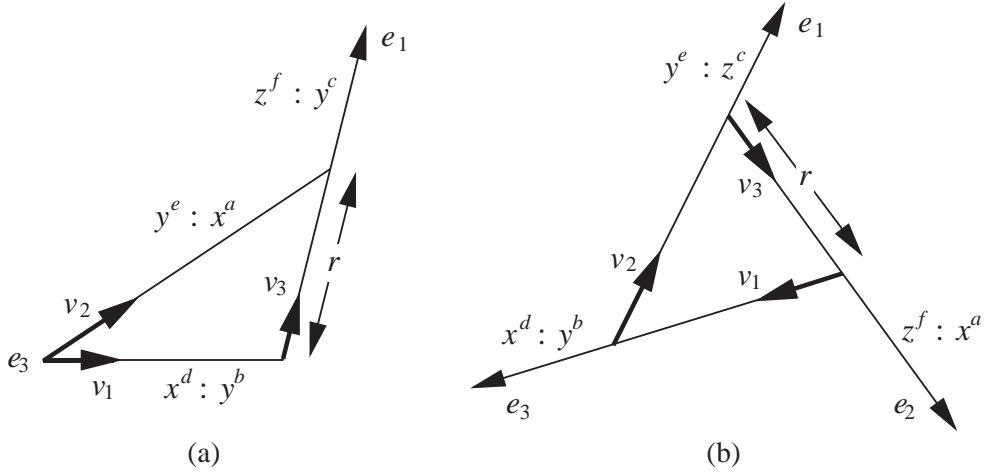}}}
 \caption{Regular triples versus monomials: (a) corner triangle; (b)
meeting of champions}
 \label{fig:2reg}
 \end{figure}

 The overlattice $L$ is based by $e_i,v_1,v_2$ for
any $i=1,2$ or 3 and any regular triple $v_1,v_2,v_3$ (or more generally
by any point of $\Z^2_\De$, together with any basis $v_1,v_2$ of the
translation lattice $\Z^2$ of $\Z^2_\De$). In contrast, $e_1,e_2,e_3$ base
$\Z^3\subset L$, and $x,y,z$ base the dual lattice $\Z^3$ of monomials on
$\C^3$. The invariant monomials form the sublattice $M\subset\Z^3$ on
which $L$ is integral, so that $M=\Hom(L,\Z)$.  Write $R$ for one of
the regular triangles of Figure~\ref{fig:2reg}. Each side of $R$
defines a sublattice (say) $\{e_3,v_1\}^\perp\cap M\iso\Z$.  The ratio
$x^d:y^b$ in Figure~\ref{fig:2reg}, or the monomial $\xi=x^d/y^b$, is
the basis of $\{e_3,v_1\}^\perp\cap M$ on which the triangle is positive, say $v_2(\xi)>0$.  (Explicit calculations are carried out for \(\frac{1}{11}(1,2,8)\) on p.~\pageref{ratiosexample}.)
 
 \begin{prop}\label{prop:1}
 Every regular triangle of side $r$ gives rise to the invariant ratios
of Figure~\ref{fig:2reg} (we permute $x,y,z$ if necessary). Moreover,
 \begin{align}
 & d-a=e-b-c=f=r \quad \text{in Case a,}
 \label{eq:1a} \\
 & d-a=e-b=f-c=r \quad \text{in Case b.}
 \label{eq:1b}
 \end{align}
\rm Note: $b,d$ (etc.)\ are not necessarily coprime; but $x^d/y^b$ is
primitive in $M$, that is, not a power of an {\em invariant} monomial.
 \end{prop}

 \begin{prop}\label{prop:2}
 Let\/ $l$ be any lattice line of\/ $\Z^2_\De$, and\/ $\bm\in M$ an
invariant monomial that bases its orthogonal\/ $l^\perp\cap M$ (as
explained at the start of Section~\ref{ssec:regtrip}). Then the
lattice lines of\/ $\Z^2_\De$ parallel to $l$ are orthogonal to
$\bm(xyz)^i$ for $i\in\Z$.

 It follows that the regular tesselations of the regular triangles of
Figure~\ref{fig:2reg} are cut out by the ratios
 \begin{align}
 & x^{d-i}:y^{b+i}z^i, \quad y^{e-j}:z^jx^{a+j}, \quad
 z^{f-k}:x^ky^{c+k} \quad \text{in Case a,}
 \label{eq:2a} \\
 & x^{d-i}:y^{b+i}z^i, \quad y^{e-j}:z^{c+j}x^j, \quad
 z^{f-k}:x^{a+k}y^k \quad \text{in Case b,}
 \label{eq:2b}
 \end{align}
for $i,j,k=0,\dots,r-1$.
 \end{prop}

 \begin{pfof}{Propositions 3.1 and 3.2} For the equalities (\ref{eq:1a}) in Case~a, note that Figure~\ref{fig:2reg}.a gives $v_1,v_2,v_3$ up to proportionality:
 \begin{equation}
 \renewcommand{\arraystretch}{1.3}
 \begin{array}{cclcr}
 v_1&\sim&(b, d, -(b+d)), \\
 v_2&\sim&(e, a, -(a+e)),\\
 v_3&\sim&(c+f, -f, -c).
 \end{array}
 \label{eq:prop}
 \end{equation}
We claim that {\em the constants of proportionality are all equal, and
equal to}
 \[
 \frac{1}{de-ab}\ =\ \frac{1}{ac+af+ef}\ =\ \frac{1}{bf+cd+df}\,.
 \]
(The denominators are the $2\times2$ minors in the array of
(\ref{eq:prop}).) For this, write
 \[
 \xi=\frac{x^d}{y^b}, \quad \eta=\frac{y^e}{x^a}, \quad
 \ze=\frac{z^f}{y^c}\,.
 \]
These 3 monomials are {\em not} a basis of $M$ (unless $r=1$, when our
regular triangle is basic). But any two of them are {\em part of a basis}.
Indeed, let $e$ be any vertex of $R$ and $\pm v_i,\pm v_j$ primitive
vectors along its two sides; then $\{e,\pm v_i,\pm v_j\}$ is a basis of
$L$, and the two monomials along the sides are part of the dual basis of
$M$. Now there are lots of dual bases around, and the claim follows at
once from
 \[
 v_1(\eta)=v_2(\xi)=v_3(\xi)=1, \quad
 v_1(\ze)=v_2(\ze)=v_3(\eta)=-1.
 \]
(The signs can be read from Figure~\ref{fig:2reg}.)  

 Equating components of $v_1+v_3=v_2$ gives $e=b+c+f$ and $a=d-f$, the
first two equalities of (\ref{eq:1a}). For the final equality, if we start
from $e_3$ and take $f$ steps along the vector $v_1$, we arrive at
 \[
 e_3+fv_1=\frac{1}{de-ab}\Bigl(bf,df,de-ab-bf-df\Bigr).
 \]
The final entry $de-ab-bf-df$ evaluates to $cd$. Thus this point has last
two entries $df,cd$ proportional to $f,c$, so lies on the third side of
$R$. Therefore $r=f$.

The proof of (\ref{eq:1b}) in Case~b is similar, and left for your
amusement. For Proposition~\ref{prop:2}, write $m,u\in M_\R$ for the
linear forms on $L$ corresponding to the monomials $\bm,xyz\in M$. The
junior plane $\R^2_\De$ is defined by $u=1$; therefore
$\{(m+iu)^\perp\}_{i\in\R}$ is a pencil of parallel lines in $\R^2_\De$.
For any lattice point $P\in\Z^2_\De$ we have $m(P)\in\Z$ and $u(P)=1$, so
$(m+iu)^\perp$ can only contain a lattice point for $i\in\Z$.
 \qed
 \end{pfof}

 \paragraph{Remark} The coordinates of points of the tesselation can be
calculated in many ways: for example, in Case~a, we get
 \[
 e_3+iv_1+jv_2
 =\frac{1}{de-ab}\bigl(bj+ei,dj+ai,de-ab-(a+e)i-(b+d)j\bigr),
 \]
which could be used to prove Proposition~\ref{prop:2}; or from the
$2\times2$ minors of
 \[
 \begin{pmatrix}
 d-i & -(b+i) & -i \\
 -(a+j) & e-j & j
 \end{pmatrix}.
 \]
It is curious that these explicit calculations in the general case shed
almost no light on Propositions~\ref{prop:1}--\ref{prop:2}, even when you
know the answers. In contrast, practice with a few numerical examples
shows at once what's going on.

\begin{figure}[!ht]
 \centering
 \input{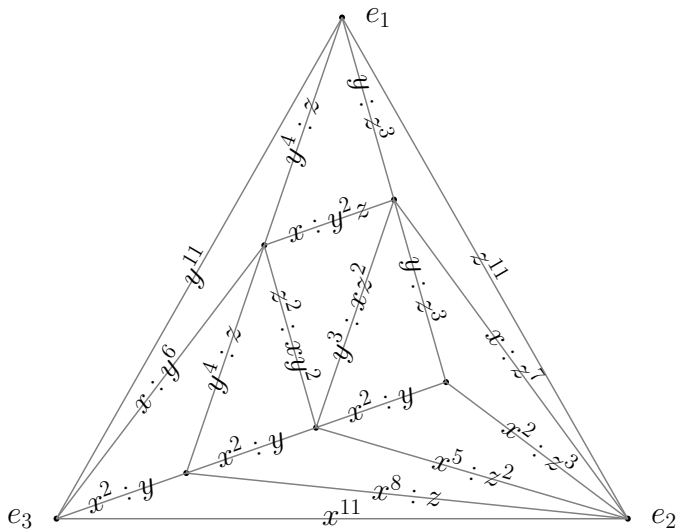}
 \caption{Ratios on the exceptional curves in $\ahilb$ for $\frac{1}{11}(1,2,8)$}
 \label{fig:ratios}
 \end{figure}

 \paragraph{Example}
 \label{ratiosexample}
 Consider once again \(\textstyle{\frac{1}{11}(1,2,8)}\).  The line from \(e_{3}\) to the lattice point
 \(\textstyle{\frac{1}{11}(1,2,8)}\) represents a 2-dimensional cone
 \(\tau\) in \(\R^{3}\) with normal vector \(\pm(2,-1,0)\).  The
 corresponding toric stratum is \(\mathbb{P}^{1}\) obtained by gluing
 \(\Spec\ \C[x^{2}y^{-1}]\) to \(\Spec\ \C[x^{-2}y]\),  so is
 parametrised by the \(A\)-invariant ratio \(x^{2}\!:\!y\).  Repeat
 for all lines to produce Figure~\ref{fig:ratios}.  

 The edges of \(\Sigma\) are not cut out by ratios;  rather,  the edges determine a single copy of \(\C\) with coordinate an invariant monomial.  That is,  the image of the \(x\),  \(y\) or \(z\)-axis of \(\C^{3}\) under the quotient map \(\pi\colon \C^{3}\to \C^{3}/A\); in this case the invariant monomials are \(x^{11}, y^{11}, z^{11}\).

 \subsection{Basic triangles and their dual monomial bases}
 \label{ssec:btdmb}
\begin{figure}[thb]
 \centerline{\mbox{\epsfbox{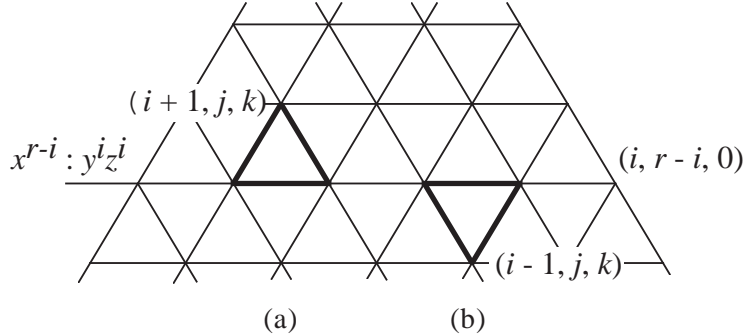}}}
 \caption{(a) Up triangle; (b) down triangle (same $i$, nonspecific $j,k$)}
 \label{fig:rr}
 \end{figure}

 The regular tesselation of a regular triangle $R$ of side $r$ is a simple
and familiar object. A moment's thought shows that every basic triangle $T$
is one of the following two types (see Figure~\ref{fig:rr} for the subgroup $\Z/r^2\subset\SL(3,\Z)$):
 \begin{description}
 \item[``up''] For $i,j,k\ge0$ with $i+j+k=r-1$, push the three sides of
$R$ inwards by $i$, $j$ and $k$ lattice steps respectively. (There are
$\binom{r+1}{2}$ choices.) We visualise three shutters closing in until
they leave a single basic triangle $T$. Note that $T$ is a scaled down
copy of $R$, parallel to $R$ and in the same orientation; in other words,
up to a translation, it is $\frac{1}{r}R$.

 \item[``down''] For $i,j,k>0$ with $i+j+k=r+1$, push the three sides of
$R$ inwards by $i$, $j$ and $k$ lattice steps (giving $\binom{r}{2}$
choices). Now the shutters close over completely, until they have a triple
overlap consisting of a single basic triangle $T$, in the opposite
orientation to $R$; up to translation, it is $-\frac{1}{r}R$.
 \end{description}

 A basic triangle $T$ has a basic dual cone in the lattice $M$, based by 3
monomials perpendicular to the 3 sides of $T$. These monomials are given
by Proposition~\ref{prop:2}, or more explicitly as follows.

 \begin{cor}\label{cor:xi}
 Let $R$ be one of the regular triangle of Figure~\ref{fig:2reg}. Its up
basic triangles have dual bases
 \begin{align*}
 & \xi=x^{d-i}/y^{b+i}z^i, \quad \eta=y^{e-j}/z^jx^{a+j}, \quad
 \ze=z^{f-k}/x^ky^{c+k} \quad\text{in Case~a} \\
 & \xi=x^{d-i}/y^{b+i}z^i, \quad \eta=y^{e-j}/z^{c+j}x^j, \quad
 \ze=z^{f-k}/x^{a+k}y^k \quad\text{in Case~b}
 \end{align*}
for $i,j,k\ge0$ with $i+j+k=r-1$. Its down basic triangles have dual bases
 \begin{align*}
 & \la=y^{b+i}z^i/x^{d-i}, \quad \mu=z^jx^{a+j}/y^{e-j}, \quad
 \nu=x^ky^{c+k}/z^{f-k} \quad\text{in Case~a} \\
 & \la=y^{b+i}z^i/x^{d-i}, \quad \mu=z^{c+j}x^j/y^{e-j}, \quad
 \nu=x^{a+k}y^k/z^{f-k} \quad\text{in Case~b}
 \end{align*}
for $i,j,k>0$ with $i+j+k=r+1$.
 \end{cor}

 \paragraph{Example $A=\Z/r\oplus\Z/r$} 
 The lattice is
 \[
 \Z^3+\Z\cdot\frac{1}{r}(1,-1,0)+\Z\cdot\frac{1}{r}(0,1,-1),
 \]
and $\De$ is spanned as usual by $e_1=(1,0,0)$, $e_2=(0,1,0)$,
$e_3=(0,0,1)$. We omit denominators as usual, writing lattice points of
$\De$ as $(a,b,c)$ with $a+b+c=r$.

 An up triangle $T$ has vertexes $(i+1,j,k)$, $(i,j+1,k)$ and $(i,j,k+1)$
for some $i,j,k\ge0$ with $i+j+k=r-1$ as in Figure~\ref{fig:rr}.a. Since
$T$ is basic, so is its dual cone in the lattice of monomials, so the dual
cone has the basis
 \[
 \xi=x^{r-i}/y^iz^i, \quad \eta=y^{r-j}/x^jz^j, \quad \ze=z^{r-k}/x^ky^k.
 \]
Thus the affine piece $Y_T=\C^3_{\xi,\eta,\ze}\subset Y_\Si$ parametrises
equations of the form
 \begin{equation}
 \begin{aligned}
 x^{r-i}&=\xi y^iz^i, \\
 y^{r-j}&=\eta x^jz^j, \\
 z^{r-k}&=\ze x^ky^k,
 \end{aligned} \qquad
 \begin{aligned}
 y^{i+1}z^{i+1}&=\eta\ze x^{r-i-1}, \\
 x^{j+1}z^{j+1}&=\xi\ze y^{r-j-1}, \\
 x^{k+1}y^{k+1}&=\xi\eta z^{r-k-1},
 \end{aligned} \qquad
 xyz=\xi\eta\ze.
 \end{equation}

 A down triangle $T$ has vertexes $(i-1,j,k)$, $(i,j-1,k)$ and $(i,j,k-1)$
for some $i,j,k\ge0$ with $i+j+k=r+1$ as in Figure~\ref{fig:rr}.b. The
sides of $T$ again correspond to the invariant ratios $x^{r-i}:y^iz^i$
etc., and its dual has basis
 \[
 \la=y^iz^i/x^{r-i}, \quad \mu=x^jz^j/y^{r-j}, \quad \nu=x^ky^k/z^{r-k}.
 \]
The affine piece $Y_T=\C^3_{\la,\mu,\nu}\subset Y_\Si$ parametrises the
equations
 \begin{equation}
 \begin{aligned}
 y^iz^i&=\la x^{r-i}, \\
 x^jz^j&=\mu y^{r-j}, \\
 x^ky^k&=\nu z^{r-k},
 \end{aligned} \qquad
 \begin{aligned}
 x^{r-i+1}&=\mu\nu y^{i-1}z^{i-1}, \\
 y^{r-j+1}&=\la\nu x^{j-1}y^{j-1}, \\
 z^{r-k+1}&=\la\mu x^{k-1}y^{k-1},
 \end{aligned} \qquad
 xyz=\la\mu\nu.
 \end{equation}

 \paragraph{Example:  regular corner triangle of side $r=1$}
The invariant ratios corresponding to the sides of a corner triangle $T$
are shown in Figure~\ref{fig:2reg}.a, where the integers
$r,a,b,c,d,e,f$ are related as in Proposition~\ref{prop:2}. If $T$ has side
$r=1$, it is basic, as is the dual cone in the lattice of monomials. The
basis consists of the invariant ratios
 \[
 \xi=x^{a+1}/y^b, \quad \eta=y^{b+c+1}/x^a, \quad \ze=z/y^c.
 \]
It follows that $\C^3_T=\C^3_{\xi,\eta,\ze}\subset Y_\Si$ parametrise the
system of equations (of which several are redundant):
 \begin{equation}
 \begin{aligned}
 x^{a+1}&=\xi y^b \\
 y^{b+c+1}&=\eta x^a \\
 z&=\ze y^c
 \end{aligned}, \quad
 \begin{aligned}
 y^{b+1}z&=\eta\ze x^a \\
 x^{a+1}z&=\xi\ze y^{b+c} \\
 xy^{c+1}&=\xi\eta
 \end{aligned}, \quad
 xyz=\xi\eta\ze.
 \label{eq:corner}
 \end{equation}

 \subsection{Remarks}
 \subsubsection{Rough proof of Theorem~\ref{th:2}}\label{sssec:rf}
 The standard construction of toric geometry is that $Y_\Si$ is the union
of the affine pieces $Y_T=\Spec k[T^\vee\cap M]$ taken over all the
triangles $T$ making up the fan $\Si$. Corollary~\ref{cor:xi} says that
$k[T^\vee\cap M]=k[\xi,\eta,\ze]$ (respectively $k[\la,\mu,\nu]$), that is,
$Y_T\iso\C^3\subset Y_\Si$, with affine coordinates $\xi,\eta,\ze$
(respectively $\la,\mu,\nu$). On the other hand Corollary~\ref{cor:xi}
also causes $Y_T$ to parametrise systems of equations such as
 \[
 x^{d-i}=\xi y^{b+i}z^i, \quad y^{e-j}=\eta z^jx^{a+j}, \quad
 z^{f-k}=\ze x^ky^{c+k}, \quad\text{etc.}
 \]
To prove Theorem~\ref{th:2}, we show that these equations determine a certain $A$-cluster of $\C^3$, and conversely, every
$A$-cluster occurs in this way; thus $Y_T$ is naturally a parameter space
for $A$-clusters. The details are given in Section~\ref{sec:pf}.

 \subsubsection{The knock-out rule \ref{sssec:KO} in exponents}
 Suppose that two lines $L_{ij}$ from the regular subdivision intersect at
an interior point of $\De$; they necessarily come out of different
vertexes, say for clarity,  $e_1$ and $e_3$. Thus they correspond to
primitive ratios $z^f:y^c$ and $y^e:x^a$. Then
 \begin{equation}
 \begin{matrix}
 \text{a line continues beyond the crossing point if and} \\
 \text{only if it has the strictly smaller exponent of $y$.}
 \end{matrix}
 \label{eq:KO}
 \end{equation}
The proof follows from Figure~\ref{fig:2reg} and the equalities of
Proposition~\ref{prop:1}; we leave the details as an exercise.

 \section{The equations of $A$-clusters}\label{sec:Nakth}

 \subsection{Two different definitions of $\GHilb M$}\label{ssec:def}
 We start with a mild warning. The literature uses two a priori different
notions of $\GHilb$: in one we set $n=|G|$, take the Hilbert scheme
$\Hilb^nM$ of all clusters of length $n$, then the fixed locus
$(\Hilb^nM)^G$, and finally, define $\GHilb M$ as the irreducible
component containing the general $G$-orbit, so birational to $M/G$. This
is a ``dynamic'' definition: a cluster $Z$ is allowed in if it is a flat
deformation of a genuine $G$-orbit of $n$ distinct points. Thus the
dynamic $\GHilb$ is irreducible by definition, but we don't really know
what functor it represents. Also, the definition involves the Hilbert
scheme $\Hilb^nM$, which is almost always very badly singular. (This point
deserves stressing: $\Hilb^nM$ is much more singular than anything needed
for $\GHilb$. As Mukai remarks, the right way of viewing $\GHilb$ should
be as a {\em variation of GIT quotient} of $X=\C^3/G$.)

Here we use the algebraic definition: a $G$-{\em cluster} $Z$ is a
$G$-invariant subscheme $Z\subset M$ with $\Oh_Z$ the regular
representation of $G$. The $G$-{\em Hilbert scheme} $\GHilb M$ is the
moduli space of $G$-clusters. Ito and Nakamura prove by continuity that a
dynamic $G$-cluster satisfies this condition, so that the dynamic
$G$-Hilbert scheme is contained in the algebraic, but the converse is not
obvious: a priori, $\GHilb M$ may have exuberant components (and quite
possibly does in general in higher dimensions). 

 Ito and Nakajima~\cite[\S2.1]{IN} prove that the algebraic and
the dynamic definitions of \(\ahilb\) coincide for a finite Abelian
subgroup \(A\subset \SL(3,\C)\).  More recently,  Bridgeland,  King
and Reid~\cite{BKR} prove that the definitions coincide
for a finite (not necessarily Abelian) subgroup \(G\subset
\SL(3,\C)\).

 \subsection{Nakamura's theorem}
 \begin{theorem}[\cite{N}]\label{th:3}
 (I) For every finite diagonal subgroup $A\subset\SL(3,\C)$ and every
$A$-cluster $Z$, generators of the ideal $\sI_Z$ can chosen as the system
of\/ $7$ equations
 \begin{equation}
 \begin{aligned}
 x^{l+1}&=\xi y^bz^f, \\
 y^{m+1}&=\eta z^cx^d, \\
 z^{n+1}&=\ze x^ay^e,
 \end{aligned} \qquad
 \begin{aligned}
 y^{b+1}z^{f+1}&=\la x^l, \\
 z^{c+1}x^{d+1}&=\mu y^m, \\
 x^{a+1}y^{e+1}&=\nu z^n,
 \end{aligned} \qquad
 xyz=\pi.
 \label{eq:Aclus}
 \end{equation}
Here $a,b,c,d,e,f,l,m,n\ge0$ are integers, and
$\xi,\eta,\ze,\la,\mu,\nu,\pi\in\C$ are constants satisfying
 \begin{equation}
 \la\xi=\mu\eta=\nu\ze=\pi.
 \end{equation}

 (II) Moreover, exactly one of the following cases holds:
 \begin{align}
 \hbox{\rm ``up''}\quad
 &\begin{cases}
 \la=\eta\ze, \quad \mu=\ze\xi, \quad \nu=\xi\eta, \quad
 \pi=\xi\eta\ze
\\
 l=a+d, \quad m=b+e, \quad n=c+f; \quad\text{or}
 \end{cases}
 \label{eq:up}
 \\[6pt]
 \hbox{\rm ``down''}\quad
 &\begin{cases}
 \xi=\mu\nu, \quad \eta=\nu\la, \quad
 \ze=\la\mu, \quad \pi=\la\mu\nu \\
 l=a+d+1, \quad m=b+e+1, \quad n=c+f+1.
 \end{cases}
 \label{eq:down}
 \end{align}
 \end{theorem}

 \paragraph{Remarks}
 The group $A$ doesn't really come into our arguments, which deal with
{\em all\/} diagonal groups at one and the same time. For example, $A=0$
makes perfectly good sense. The particular group for which $Z$ is an
$A$-cluster is determined from the exponents in (\ref{eq:Aclus}) as
follows: its character group $A^*$ is generated by its eigenvalues
$\chi_x,\chi_y,\chi_z$ on $x,y,z$, and related by
 \begin{equation}
 \chi_x+\chi_y+\chi_z=0 \quad\text{and}\quad\ 
 \begin{aligned}
 (l+1)\chi_x&=b\chi_y+f\chi_z \\
 (m+1)\chi_y&=c\chi_z+d\chi_x \\
 (n+1)\chi_z&=a\chi_x+e\chi_y.
 \end{aligned}
 \label{eq:char}
 \end{equation}
This is a presentation of $A$ as a $\Z$-module, as a little $4\times3$
matrix; all our stuff about regular triples, regular tesselations and so
on, can be viewed as a classification of different presentations of $A^*$
of type (\ref{eq:char}).

 The equations of $Z$ in Theorem~\ref{th:3} may be redundant (for example,
(\ref{eq:corner})), and the choice of exponents $a,b,\dots,n$ is usually
not unique: a cluster with $\pi\ne0$ corresponds to a point in the big
torus of $Y_\Si$, belonging to every affine set $Y_T$, and thus can be
written in {\em every} form consistent with the group $A$.

 Although at this point we're sober characters doing straight-laced
algebra, the argument is substantially the same as that already sketched in
\cite{R}, which you may consult for additional examples, pictures,
philosophy and jokes. See also \cite{N}.

 \begin{pfof}{(I)} By definition (see \ref{ssec:def}), the Artinian ring
$\Oh_Z=k[x,y,z]/I_Z=\Oh_{\C^3}/\sI_Z$ of $Z$ is the regular
representation, so each character of $A$ has exactly a one dimensional
eigenspace in $\Oh_Z$. Arguing on the identity character and using the
assumption $A\subset\SL(3,\C)$ provides an equation $xyz=\pi$ for some
$\pi\in\C$.

 Since $k[x,y,z]$ is based by monomials, their images span $\Oh_Z$;
monomials are eigenfunctions of the $A$ action. Obviously, each eigenspace
in $\Oh_Z$ contains a nonzero image of a monomial $\bm$, and is based by
any such. Moreover, if $\bm$ is a multiple of an invariant monomial, say
$\bm=\bm_0\bm_1$ with $\bm_0$ invariant under $A$, and is nonzero in
$\Oh_Z$, then the other factor $\bm_1$ is also a basis of the same
eigenspace. {From} now on, we say {\em basic monomial\/} in $\Oh_Z$ to mean
the nonzero image in $\Oh_Z$ of a monomial that is not a multiple of an
invariant monomial; in particular, it is not a multiple of $xyz$, so
involves at most two of $x,y,z$.

 The next result shows how to choose the equations in (\ref{eq:Aclus}). 

 \begin{lemma} Let $x^r$ be the first power of $x$ that is $A$-invariant.
Then there is (at least) one $l\in[0,r-1]$ such that\/
$1,x,x^2,\dots,x^l\in\Oh_Z$ are basic monomials, and\/ $x^{l+1}$ is a
multiple of some basic monomial $y^bz^f$ in the same eigenspace, say
$x^{l+1}=\xi y^bz^f$ for some $\xi\in\C$.
 \end{lemma}

Let's first see that the lemma gives the equations in (I). Indeed
$x^{l+1},y^bz^f$ belong to a common eigenspace, and therefore, because
$xyz$ is invariant, also $x^l$ and $y^{b+1}z^{f+1}$ belong to a common
eigenspace. This is based by $x^l$ by choice of $l$, hence we get the
relation $y^{b+1}z^{f+1}=\la x^l$.

 Finally, since $y^bz^f$ is a basic monomial, $\la\xi=\pi$ corresponds to
the syzygy $\la(\one)+x(\two)-y^bz^f(\three)$ between the three relations
 \[
 (\one)\quad x^{l+1}=\xi y^bz^f, \quad
 (\two)\quad y^{b+1}z^{f+1}=\la x^l, \quad
 (\three)\quad xyz=\pi.
 \]
The relations involving $y^{m+1}$ and $z^{n+1}$ arise similarly.

 \begin{pfof}{the lemma} If $x^{r-1}\ne0\in\Oh_Z$ it is a basic monomial,
and one choice is to take $l=r-1$ and $b=f=0$, and to take the relation
$x^{l+1}=x^r=\xi\cdot1$. (Other choices arise if the eigenspace of some
$x^{l'+1}$ with $l'<l$ also contain a basic monomial $y^{b'}z^{f'}$.)

 If not, there is some $l$ with $0\le l\le r-1$ such that
$1,x,x^2,\dots,x^l$ are basic monomials and $x^{l+1}=0\in\Oh_Z$. Now the
eigenspace of $x^{l+1}$ must contain a basic monomial $\bm$; under the
current assumptions, we assert that $\bm$ is of the form $y^bz^f$, which
proves the lemma. We need only prove that $\bm$ is not a multiple of $x$.
If $\bm=x\bm'$ then $\bm'$ must in turn be a basic monomial in the same
eigenspace as $x^l$. But then $x^l=(\text{unit})\cdot\bm'$ contradicts
$x^{l+1}=0$ and $x\bm'\ne0$. \qed \end{pfof}

Now (I) says that, for any $A$ and any $A$-cluster $Z$, once the relations
(\ref{eq:Aclus}) are derived as above, $\Oh_Z$ is based by the monomials
in the tripod of Figure~\ref{fig:tripod}, and the relations reduce any
monomial $\bm$ to one of these. We derived the relations in pairs
$x^{l+1}\mapsto y^bz^f$ and $y^{b+1}z^{f+1}\mapsto x^l$. The first type
 \begin{figure}[thb]
 \centerline{\mbox{\epsfbox{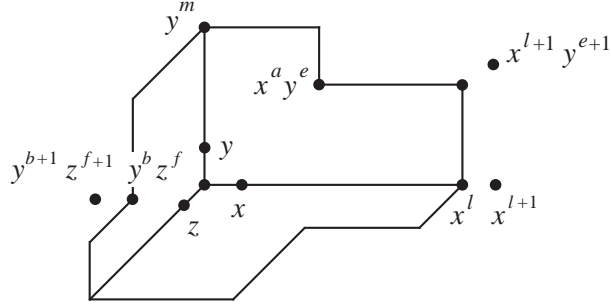}}}
 \caption{Tripod of monomials basing $\Oh_Z$}
 \label{fig:tripod}
 \end{figure}
reduces pure powers of $x$ higher than $x^l$. Suppose we have a further
relation in the first quadrant, (say) $x^\al y^\ep\mapsto \bm$: if $\bm$
involves $x$ or $y$ the new relation would be a multiple of a simpler
relation. On the other hand, if $\bm=z^\ga$ is a pure power of $z$, the
above argument shows the new relation is paired with a relation
$z^{\ga+1}\mapsto x^{\al-1} y^{\ep-1}$, which contradicts our choice of
$n$ (in the exponent of $z^{n+1}$). This concludes the proof of (I). \qed
\end{pfof}

\begin{pfof}{(II)} The point is that a monomial just off one of the
shoulders of the tripod of Figure~\ref{fig:tripod} such as $x^{l+1}y^{e+1}$
or $y^{m+1}z^{f+1}$, etc., reduces to a basic monomial in two steps
involving two of the $\xi,\eta,\ze$ relations, or two of the $\la,\mu,\nu$
relations. (Compare \cite{R}, Remark~7.3 for a discussion.)

The first reduction applies if $b+e\ge m$:
 \[
 x^{l+1}y^{e+1} \mapsto \xi y^{b+e+1}z^f
 \mapsto \xi\eta y^{b+e-m} x^dz^{c+f}
 \]
This implies that the monomials $x^{l-d+1}y^{m-b+1}$ and $z^{c+f}$ are in
the same eigenspace, and the existence of the relation
 \[
 x^{l-d+1}y^{m-b+1}=\xi\eta z^{c+f}
 \]
between them. But from the argument in (I), there is only one relation in
this quadrant, namely $x^{a+1}y^{e+1}=\nu z^n$. Therefore $l-d=a$,
$m-b=e$, $c+f=n$ and $\nu=\xi\eta$. Now $a+d\ge l$ and $c+f\ge n$, so that
we can run the same two-step reduction to other monomials to get
$\la=\eta\ze$ and $\mu=\xi\ze$.

The second type of reduction applies if $m\ge b+e+1$
 \[
 y^{m+1}z^{f+1} \mapsto \la y^{m-b}x^l
 \mapsto \la\nu x^{l-a-1}y^{m-b-e-1} z^n
 \]
Therefore the two monomials $y^{b+e+2}$ and $x^{l-a-1}z^{n-f-1}$ are in
the same eigen\-space, and $y^{b+e+2}=\la\nu x^{l-a-1}z^{n-f-1}$. As
before, this must be identical to the $\eta$ relation, so that $m+1=b+e+2$,
$l-a-1=d$, $n-f-1=c$ and $\eta=\la\nu$. This proves the theorem.
 \qed \end{pfof}

 \section{Proof of Theorem~\ref{th:2}}\label{sec:pf}
 The point is to identify the objects in the conclusion of
Corollary~\ref{cor:xi} and of Theorem~\ref{th:3}; this is really just a
mechanical translation. To distinguish between the two sets of symbols, in
the monomial bases of Corollary~\ref{cor:xi}, we first substitute for
$d,e,f$ from (\ref{eq:1a}--\ref{eq:1b}) of Proposition~\ref{prop:1}, and
then replace
 \[
 a\mapsto A, \quad b\mapsto B, \quad c\mapsto C.
 \]
Each of the monomial bases of Corollary~\ref{cor:xi} gives rise to a triple
of equations, either up:
 \begin{align*}
 & x^{A+r-i}=\xi y^{B+i}z^i, \quad y^{B+C+r-j}=\eta z^jx^{A+j}, \quad
 z^{r-k}=\ze x^ky^{C+k} \quad\text{in Case~a} \\
 & x^{A+r-i}=\xi y^{B+i}z^i, \quad y^{B+r-j}=\eta z^{C+j}x^j, \quad
 z^{C+r-k}=\ze x^{A+k}y^k \quad\text{in Case~b}
 \end{align*}
with $i,j,k\ge0$ and $i+j+k=r-1$; or down:
 \begin{align*}
 & y^{B+i}z^i=\la x^{A+r-i}, \quad z^jx^{A+j}=\mu y^{B+C+r-j}, \quad
 x^ky^{C+k}=\nu z^{r-k} \quad\text{in Case~a} \\
 & y^{B+i}z^i=\la x^{A+r-i}, \quad z^{C+j}x^j=\mu y^{B+r-j}, \quad
 x^{A+k}y^k=\nu z^{C+r-k} \quad\text{in Case~b}
 \end{align*}
with $i,j,k>0$ and $i+j+k=r+1$.

Each triple can be completed to the equations of an $A$-cluster; for
example, the first triple gives:
 \[
 \renewcommand{\arraycolsep}{.2em}
 \begin{array}{rcl}
 x^{A+r-i}&=&\xi y^{B+i}z^i \\
 y^{B+C+r-j}&=&\eta z^jx^{A+j} \\
 z^{r-k}&=&\ze x^ky^{C+k}
 \end{array}
 \quad
 \begin{array}{rcl}
 y^{B+r-j-k}z^{r-j-k}&=&\eta\ze x^{A+j+k} \\
 z^{r-i-k}x^{A+r-i-k}&=&\ze\xi y^{B+C+k+i} \\
 x^{r-i-j}y^{C+r-i-j}&=&\xi\eta z^{i+j}
 \end{array}
 \quad
 xyz=\xi\eta\ze.
 \]
(The method is to multiply together any two of the equations and cancel
common factors.) Since $i+j+k=r-1$, these are of the form of
Theorem~\ref{th:3}, with $l=A+j+k$, $b=B+i$,$f=i$, etc.. The other cases
are similar. Therefore as explained in \ref{sssec:rf}, each affine piece
$Y_T\iso\C^3\subset Y_\Si$ parametrises
$A$-clusters.

Conversely, we prove that for $A\subset\SL(3,\C)$ a finite diagonal
subgroup and $Z$ an $A$-cluster with equations as in Theorem~\ref{th:3},
$Z$ belongs to one of the families parametrised by $Y_T$. If $Z$ is ``up''
its equations are determined by the first three:
 \begin{equation}
 x^{a+d+1}=\xi y^bz^f,\quad y^{b+e+1}=\eta z^cx^d,\quad
 z^{c+f+1}=\ze x^ay^e.
 \label{eq:3up}
 \end{equation}
Consider first just two of the possibilities for the signs of $f-b$,
$d-c$, $e-a$.

\begin{enumerate}
 \item Suppose $b\ge f$, $d\ge c$ and $e\ge a$. We define $A,B,C,i,j,k$ by
 \[
 A=d-c, \quad B=b-f, \quad C=e-a, \quad
 i=f,\quad j=c,\quad k=a
 \]
and set $r=i+j+k+1$. Then, obviously,
 \[
 a=k, \quad b=B+i, \quad c=j, \quad d=A+j, \quad e=C+k, \quad f=i.
 \]
Substituting these values in the exponents of (\ref{eq:3up}), puts
the equations of $Z$ in the form up, Case~a.

 \item Similarly, if $b\ge f$, $c\ge d$ and $a\ge e$, we fix up
$A,B,C,i,j,k$ so that
 \[
 a=A+k, \quad b=B+i, \quad c=C+j, \quad d=j, \quad e=k, \quad f=i.
 \]
Substituting in (\ref{eq:3up}), shows that $Z$ is up, Case~b.
\end{enumerate}

One sees that the permutation $y\bij z$ leads to $b\bij f$, $a\bij d$ and
$c\bij e$, and the other possibilities for the signs of $e-a$, $f-b$,
$d-c$ all reduce to these two cases on permuting $x,y,z$. In fact,
Figure~\ref{fig:2reg}.a has 6 different images on permuting $x,y,z$
(corresponding to the choices of $e_1$ and $e_3$), and
Figure~\ref{fig:2reg}.b has 2 different images (corresponding to the
cyclic order).

If $Z$ is ``down'' its equations can be deduced from the second three:
 \begin{equation}
 y^{b+1}z^{f+1}=\la x^{a+d+1}, \quad z^{c+1}x^{d+1}=\mu x^{b+e+1},
 \quad x^{a+1}y^{e+1}=\nu z^{c+f+1}
 \label{eq:3down}
 \end{equation}
Exactly as before, if $b\ge f$, $d\ge c$ and $e\ge a$ then we can fix up
$A,B,C\ge0$ and $i,j,k>0$ so that
 \begin{gather*}
 a+1=k, \quad b+1=B+i, \quad c+1=j, \\
 d+1=A+j, \quad e+1=C+k, \quad f+1=i,
 \end{gather*}
which puts (\ref{eq:3down}) in the form down, Case~a. The rest of the proof
is a routine repetition. This proves Theorem~\ref{th:2}. \qed

 \bigskip

\noindent
Miles Reid, \\
Math Inst., Univ. of Warwick, \\
Coventry CV4 7AL, England \\
Miles@Maths.Warwick.Ac.UK

\medskip
\noindent
Alastair Craw, \\
Math Inst., Univ. of Warwick, \\
Coventry CV4 7AL, Great Britain \\
Craw@Maths.Warwick.Ac.UK

 \newpage

 \thispagestyle{empty}
 \begin{figure}[p]
 \epsfbox[40 0 450 720]{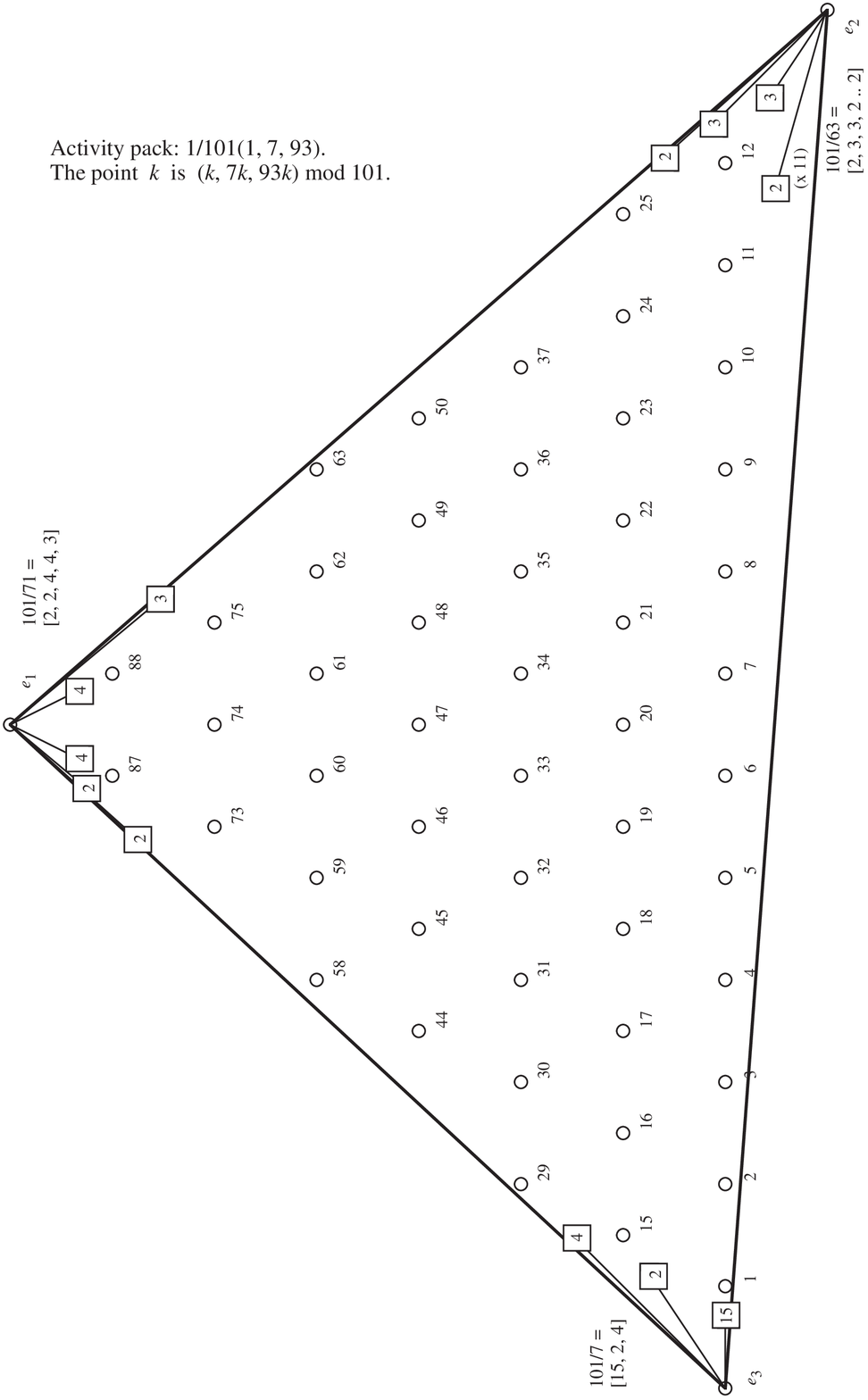}
 \label{fig:BigExa}
 \end{figure}

 \end{document}